%% file: root.tex
\documentclass[letterpaper, 10 pt, conference]{ieeeconf}  

\usepackage{amsfonts}
\usepackage{amsmath}
\usepackage{amssymb}

\usepackage[dvipsnames]{xcolor}
\usepackage{todonotes}
\usepackage{graphicx}
\usepackage{url}
\usepackage{cite}

\IEEEoverridecommandlockouts                              

\title{\LARGE \bf
Risk-Aware Control and Optimization for High-Renewable Power Grids
}

\author{Neil Barry, Minas Chatzos, Wenbo Chen, Dahye Han, Chaofan Huang, V. Roshan Joseph, Michael Klamkin,\\Seonho Park, Mathieu Tanneau, Pascal Van Hentenryck, Shangkun Wang, Hanyu Zhang and Haoruo Zhao
\thanks{Industrial and Systems Engineering, Georgia Institute of Technology, Atlanta, GA, United States.
        Corresponding author: {\tt\small mathieu.tanneau@gatech.edu}}%
\thanks{This research is partly funded by NSF Awards 2007095 and
2112533, and ARPA-E Perform Award AR0001136.}
}

\begin{document}

\maketitle
\thispagestyle{empty}
\pagestyle{empty}


\input{tex/abstract}


\input{tex/introduction}

\input{tex/forecast}

\input{tex/optimization}
\input{tex/learning}

\input{tex/risk}

\input{tex/conclusion}

\addtolength{\textheight}{-8cm}   





\bibliographystyle{IEEEtran}
\bibliography{IEEEabrv,refs}

\end{document}

%% file: tex/abstract.tex
\begin{abstract}
The transition of the electrical power grid from fossil fuels to renewable sources of energy raises fundamental challenges to the market-clearing algorithms that drive its operations. Indeed, the increased stochasticity in load and the volatility of renewable energy sources have led to significant increases in prediction errors, affecting the reliability and efficiency of existing deterministic optimization models. The RAMC project was initiated to investigate how to move from this deterministic setting into a risk-aware framework where uncertainty is quantified explicitly and incorporated in the market-clearing optimizations. Risk-aware market-clearing raises challenges on its own, primarily from a computational standpoint. This paper reviews how RAMC approaches risk-aware market clearing and presents some of its innovations in uncertainty quantification, optimization, and machine learning. Experimental results on real networks are presented.
\end{abstract}

%% file: tex/introduction.tex
\section{INTRODUCTION}
\label{sec:introduction}

\begin{figure}[!t]
    \centering
    \includegraphics[height=0.44\columnwidth]{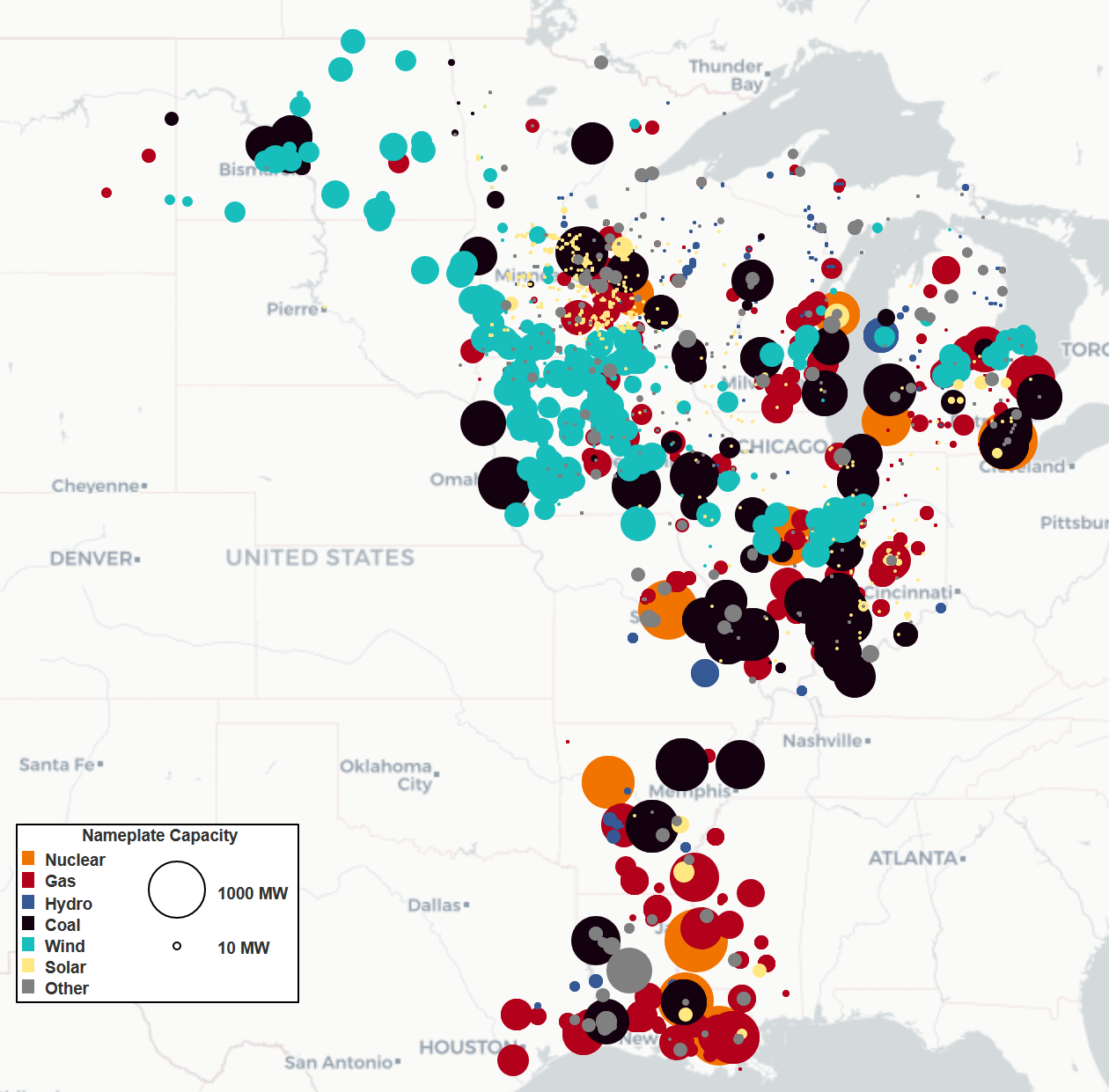}
    \hfill
    \includegraphics[height=0.44\columnwidth]{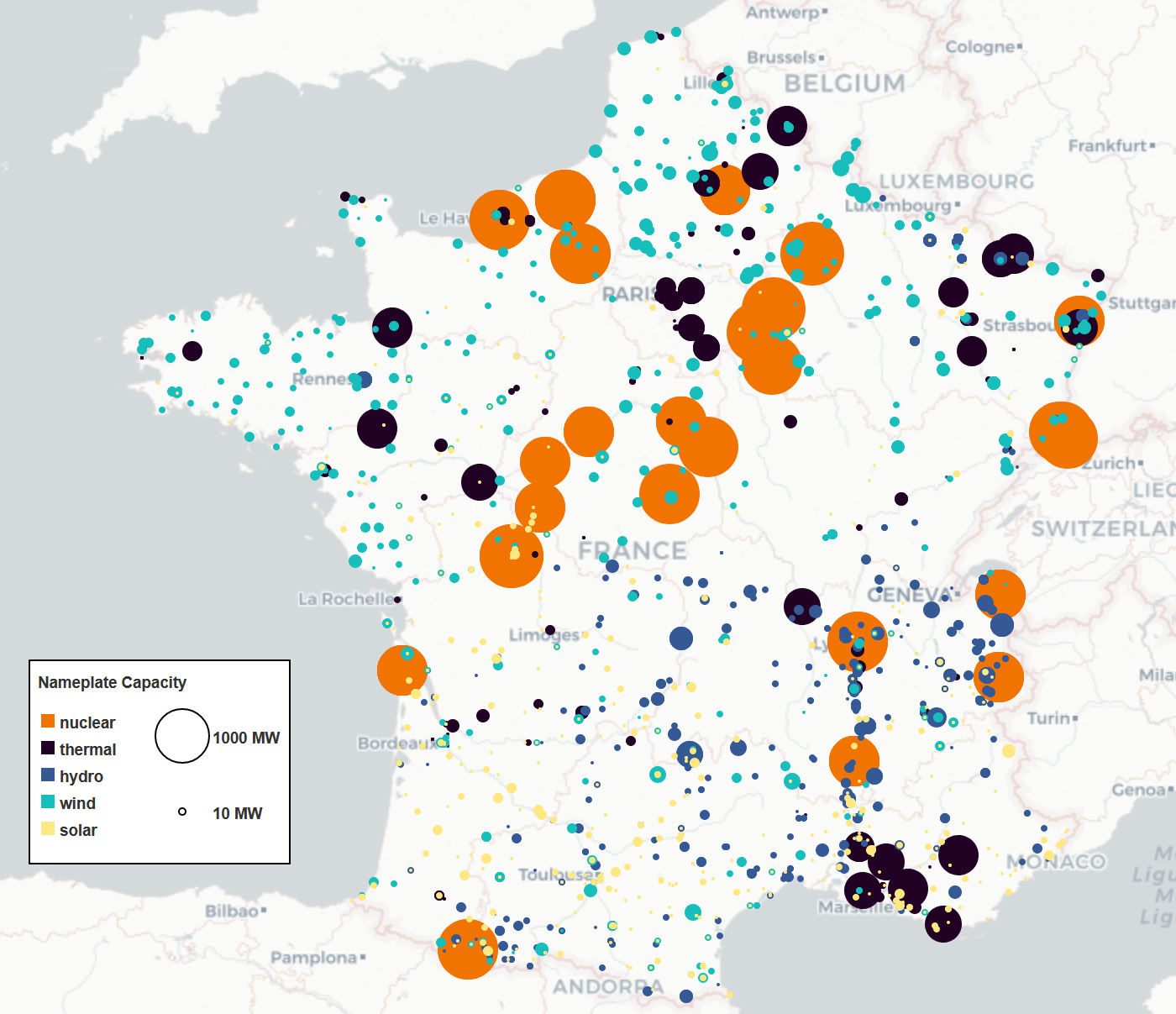}
    \caption{Size, type, and location of generation units in MISO (left, US generators only) and RTE (right) systems. Approximate locations of US and French generators obtained from \cite{EIA860} and \cite{ODRE}.}
    \label{fig:introduction:MISO_RTE_grids}
\end{figure}

Power systems are critical infrastructures for the social and economic welfare of society as a whole.
They are driven by complex optimization algorithms whose goals are to balance, at all times, the load (or demand) and generation, minimizing costs while satisfying the laws of physics and engineering constraints on generators and transmission equipment.
In many countries around the world, Independent Systems Operators (ISOs) or Transmission Systems Operators (TSOs) are responsible for the optimization algorithms driving the operations of electrical power systems.
These optimizations are challenging computationally, since they need to model nonlinear nonconvex physical constraints, as well as discrete commitment variables that decides which generators will deliver power at different times of the day.
These optimization models are also large-scale, which can be seen by the number of generators
for MISO, the Midcontinent Independent Systems Operator, and France in Figure  \ref{fig:introduction:MISO_RTE_grids}.

In the United States, the power grid is operated by a series of market-clearing algorithms at different time granularities.
They include the day-ahead unit commitment, the look-ahead commitment, and the real-time dispatching optimizations.
These market-clearing algorithms are deterministic and co-optimize energy and reserves, in order to ensure both efficiency and reliability.
The transition of the electrical power grid from fossil fuels to renewable sources of energy raises fundamental challenges to these market-clearing algorithms.
Indeed, the increased stochasticity in load and the volatility of renewable energy sources have led to significant increases in prediction errors, potentially affecting the reliability and efficiency of existing deterministic optimization models.
As the share of renewable energy in the energy mix increases, it becomes increasingly challenging to adapt existing deterministic operations to meet reliability and economic objectives. 

The RAMC project was initiated to investigate how to move from this deterministic setting into a risk-aware framework, where uncertainty is quantified explicitly and incorporated in the market-clearing optimizations.
Risk-aware market clearing raises challenges on its own obviously, primarily from a computational standpoint.
This paper reviews how RAMC approaches risk-aware market clearing and presents some of its innovations in uncertainty quantification, optimization, and machine learning (ML).
In particular, the paper shows how asset bundling and scenario generation based on support points provide a critical link between forecasting and risk-aware optimization.
It also indicates the promise and viability of stochastic optimization to hedge against risk, while minimizing costs.
Moreover, it demonstrates the potential of ML in addressing the computational challenges of risk-aware market clearing, introducing novel frameworks that include just-in-time learning, physics-informed learning based on Lagrangian duality and decomposition, and the predict-then-repair methodology.
The RAMC methodology has been validated  on real networks, and preliminary experimental are presented to illustrate the potential benefits.

The rest of the paper is structured as follows.
Section \ref{sec:forecast} reviews forecasting and uncertainty-quantification models for load and renewable output prediction.
Section \ref{sec:optimization} describes the current optimization pipeline of MISO, and discusses several extensions.
Section \ref{sec:learning} introduces novel learning methodologies.
Section \ref{sec:risk} presents further risk analysis and estimation, and Section \ref{sec:conclusion} concludes the paper.

%% file: tex/forecast.tex
\section{PREDICTIONS}
\label{sec:forecast}

\subsection{Forecasting Models}
\label{sec:predictions:forecast}

    The introduction of significant renewable sources of energy in the grid raises significant forecasting challenges. Wind and solar energy sources fluctuate according to exogenous factors, e.g., cloud coverage for solar power generation or wind speed/direction for wind power generation. Obtaining accurate forecasts for renewables is thus critical to ensure reliable and efficient operations: they are used as inputs to every market-clearing algorithms, albeit for different time horizons and discretizations. 
    
    \begin{figure}[t]
      \centering
        \includegraphics[width=0.49\columnwidth]{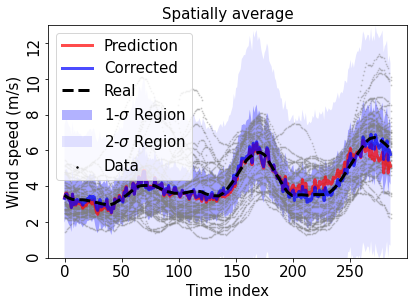}
        \includegraphics[width=0.49\columnwidth]{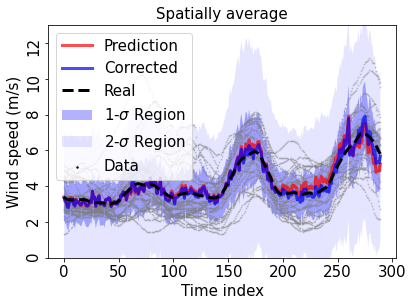}
        \caption{Examples of corrected prediction of wind speed suggested by MRSTK with two different data resolutions in MISO region. (Left: spatial resolution $\kappa=50$, temporal resolution $\tau=24$, Right: $\kappa=30$, $\tau=20$). Excerpted from \cite{Zhu2021_MultiResolutionSpatioTemporal}.}\label{fig:mrstk}
    \end{figure}
    
    The RAMC project aims at exploiting the significant spatial and temporal correlations that exists in wind and solar time series. A step in that direction was the design and implementation of a Multi Resolution Spatio-Temporal Kriging (MRSTK) model \cite{Zhu2021_MultiResolutionSpatioTemporal} for forecasting the wind speed and direction. The model relies on the dynamic graph that uses the wind direction to capture spatial and temporal correlation. The overall model consolidates individual regression models for each resolution through the Kriging model and is trained using a variational method for scalability. Its effectiveness was demonstrated on historical weather data for the MISO wind as illustrated in Figure~\ref{fig:mrstk}.
    
    The project is also evaluating, in a systematic way, the strengths and weaknesses of a variety of forecasting methods. They include various time-series models such as Periodic Vector Regression (PVAR) \cite{ursu2009modelling} and Recurrent Neural Networks (RNN) based models such as Long-Short Term Memory (LSTM) \cite{hu2018nonlinear} and Gated Recurrent Unit (GRU) \cite{chung2014empirical}.
    It is not entirely clear which forecasting models are the most accurate for the various horizons and discretizations, and how to best use the massive data sets now available.
    One intriguing direction is the use of spatio-temporal graph neural networks that integrate some atmospheric physics present in Numerical Weather Prediction (NWP) models (e.g., \cite{sweeney2020future, yu2020superposition, park2019physics}). 

\subsection{Asset Bundling}
\label{sec:predictions:asset_bundling}
    
    The RAMC project is also pioneering some novel techniques in asset bundling, i.e., in selecting which assets to predict together to facilitate ML and minimize forecasting errors.
    The benefits of asset bundling are two-fold.
    First, bundling wind farms together reduces the dimensions of the forecasting problems, thus making ML models easier to train.
    Second, bundled time series, if selected carefully, exhibit smoother patterns, which improve the accuracy of the predictions.
    The forecasting methods discussed earlier can be applied to bundled assets.
    This is then followed by a disaggregation method (e.g., \cite{chatzos2021data}) that yields forecasts at the original granularity (e.g., local balancing authority (LBA) in the MISO pipeline). 
    
    The project developed a new concept called {\em intermittency index-based bundling}.
    Its key idea is to minimize the covariance of the first-order time derivative of the bundles, while including spatial constraints to avoid bundling assets that are too far of each other.
    The optimal bundles can be obtained by solving a mixed-integer quadratic programming (MIQP) problem or approximated by a fast greedy heuristic.

\subsection{Uncertainty Quantification}
\label{sec:predictions:undertainty_quantification}
    
    Although ISOs' existing market-clearing pipelines are based on point forecasting, it is necessary to quantify uncertainty and generate scenarios for their risk-aware counterparts.
    The uncertainty quantification can inform both the energy dispatch and the reserve allocation, again ensuring both reliable and cost-effective commitments and dispatches. 
    
    MC dropout \cite{gal2016dropout} has been proposed to generate scenarios and quantify uncertainty. However, it only measures \textit{epistemic} uncertainty, i.e., the variations in the data that the model does not capture well, and in a relative manner which is not ideal in this context. What is needed is to also capture \textit{aleatoric} uncertainty, which can be obtained by the posterior distribution of a  probabilistic forecasting method. The RAMC project uses Gaussian copula model \cite{werho2021scenario}, Gaussian mixture model \cite{zhang2019short}, and normalizing flow based model \cite{dumas2022deep} for this purpose. 
        
    \subsection{Support Points}
    \label{sec:predictions:support_points}
    
    \begin{figure}[t]
      \centering
        \includegraphics[width=0.49\columnwidth]{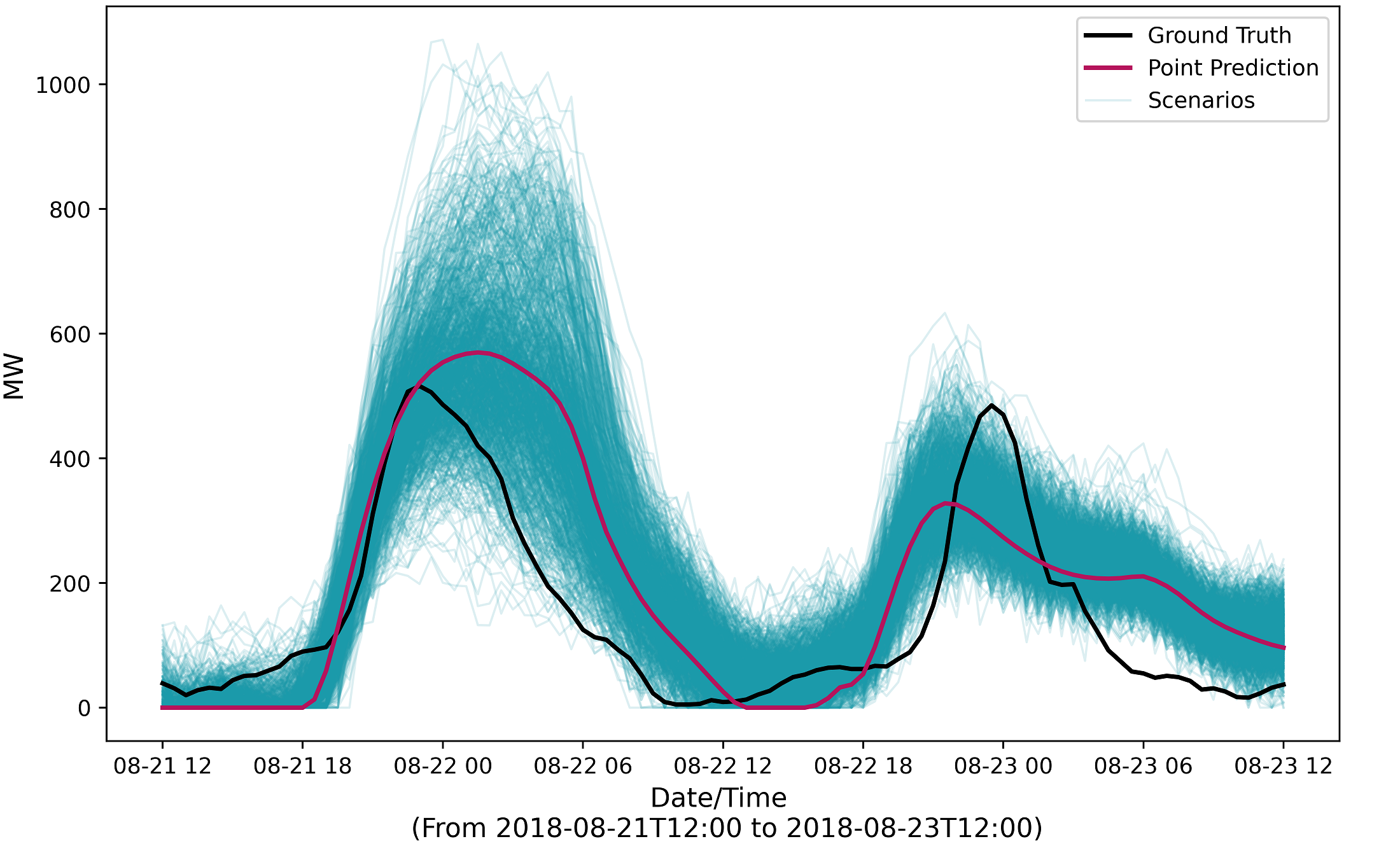}
        \includegraphics[width=0.49\columnwidth]{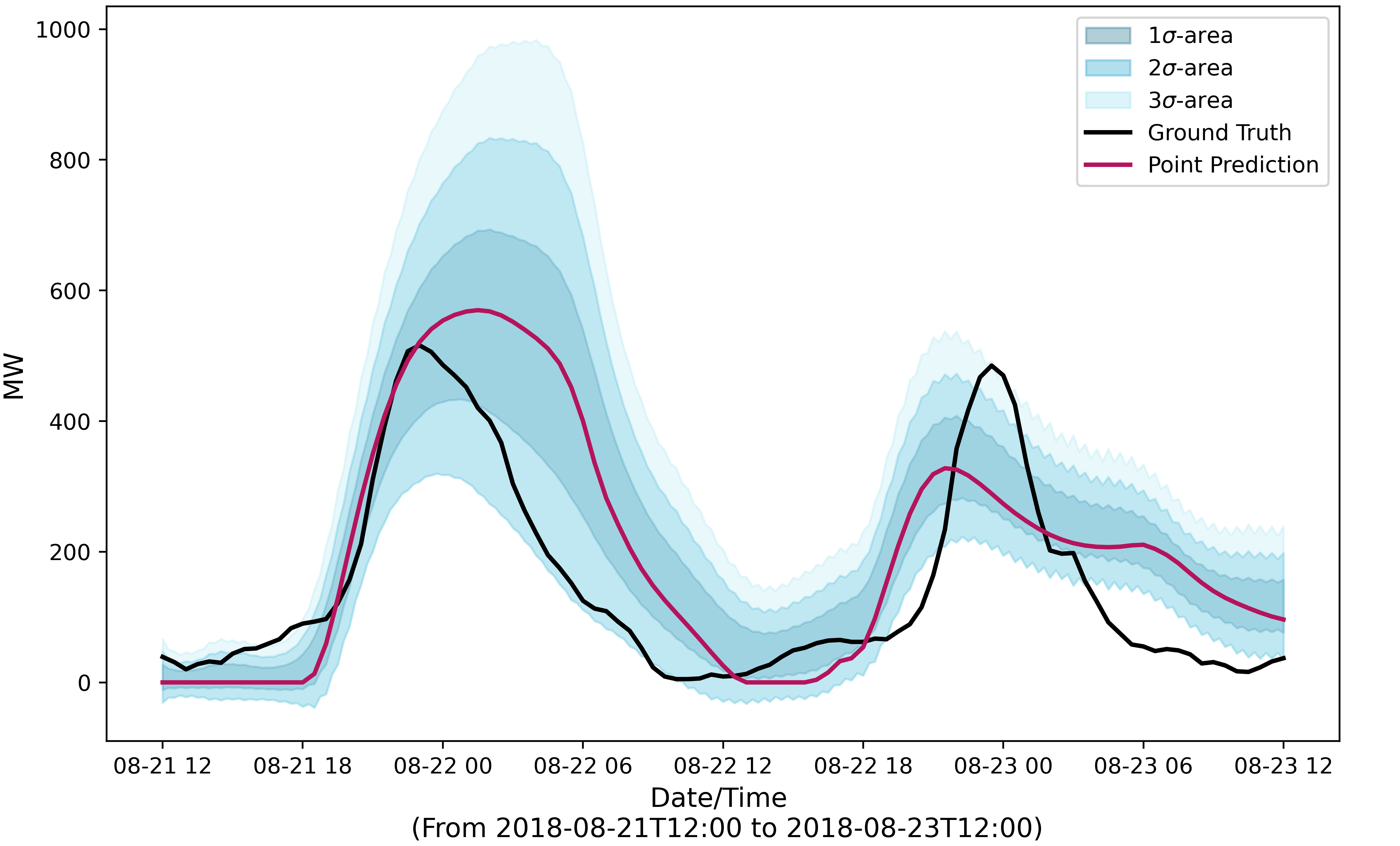}\\
        \includegraphics[width=0.49\columnwidth]{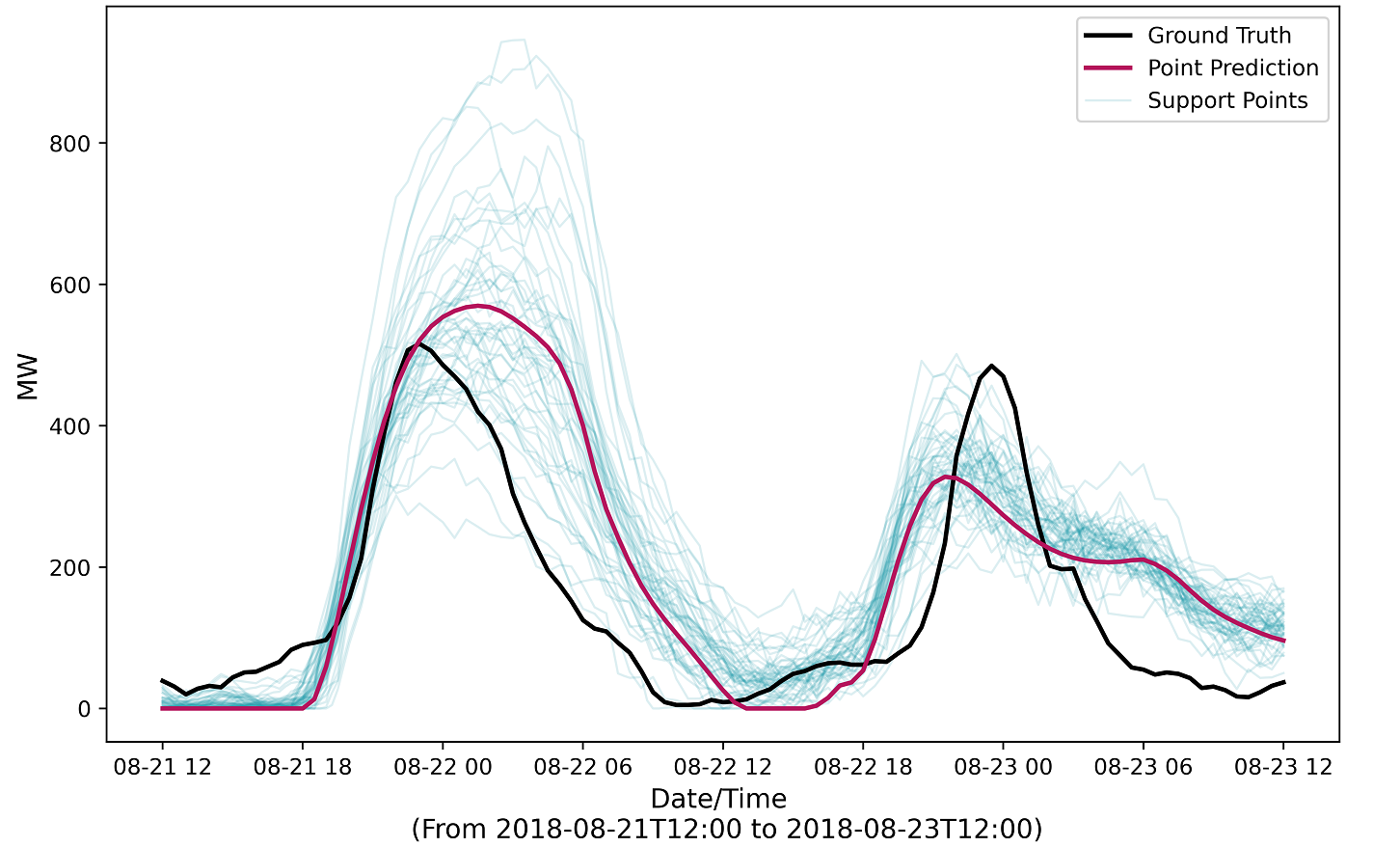}
        \includegraphics[width=0.49\columnwidth]{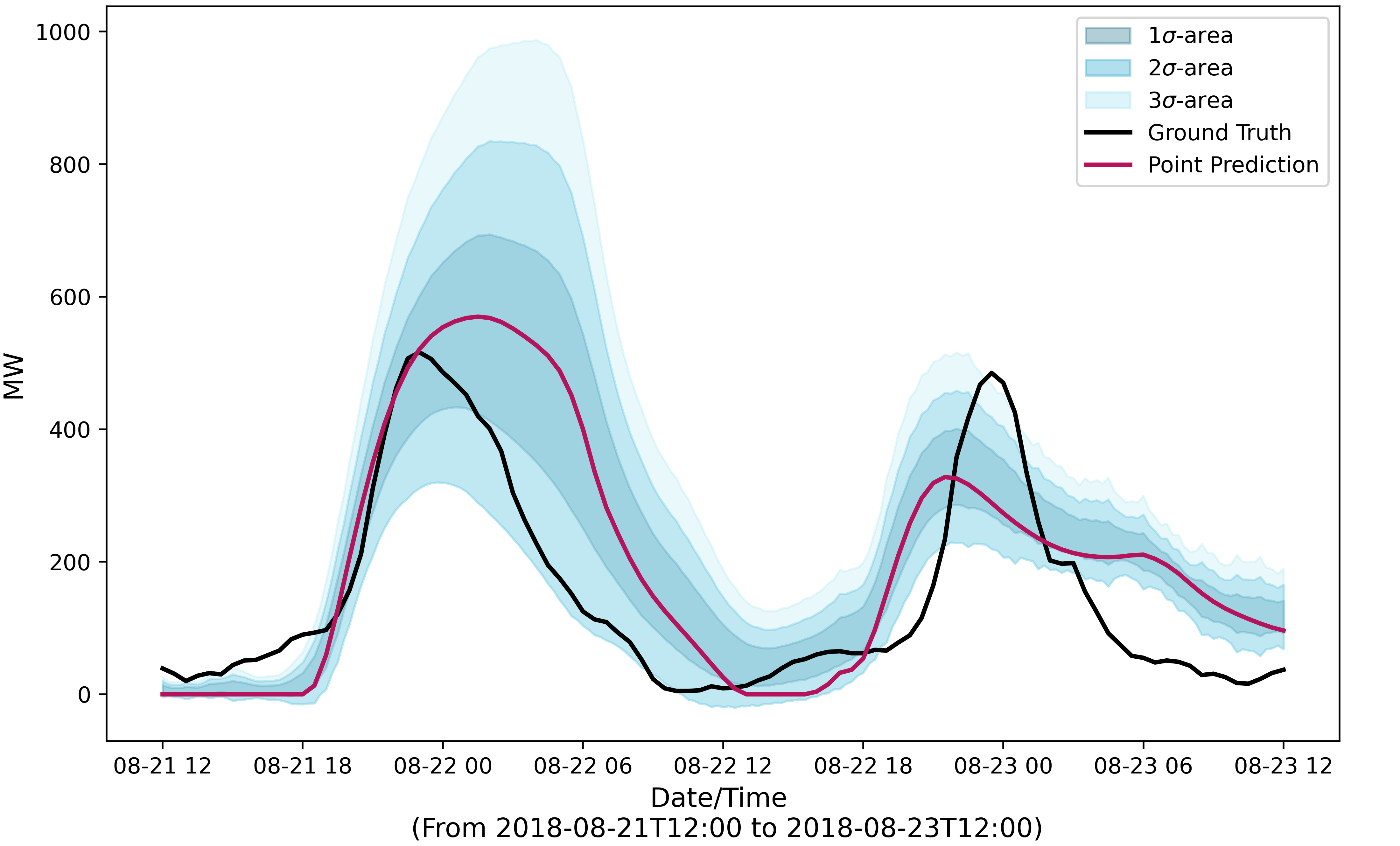}
      \caption{48 hours wind generation forecast scenarios (Top) and their support points (Bottom) on Aug. 22, 2018 in Grand-Est region in the French system. \textbf{Left}: showing the raw scenarios and support points. \textbf{Right}: showing the 1, 2, and 3 sigma levels given the scenarios and support points.}
      \label{fig:support_points}
    \end{figure}    

    The scenarios generated by the uncertainty quantification methods can be used as inputs to risk-aware market-clearing optimizations (see Section \ref{sec:optimization:stochastic} and reference \cite{zheng2014stochastic} for instance).
    The computational complexity of risk-aware optimization may be strongly affected by the number of scenarios.
    To minimize the computational requirements, the number of scenarios can be reduced by Monte-Carlo (MC) sampling, starting from a very large set of scenarios.
    The project instead leverages the concept of support points \cite{mak2018support}, which are used to select the ``best'' $k$ scenarios, i.e., scenarios that best capture the original set of scenarios while being as diverse as possible.
    Illustrative examples of support points are in Figure~\ref{fig:support_points} in which it is shown that support points capture the distribution of the forecasts with the smaller number of scenarios.
    Support points are obtained by solving a convex-concave optimization for which fast and high-fidelity heuristics are available \cite{Vakayil22_DataTwinning}.
    The project has generated support scenarios from a posterior distribution produced by Gaussian Process model \cite{chen2013wind} and a large set of scenarios produced by a Gaussian copula model \cite{werho2021scenario}.
    The comparison between two approaches has been conducted by measuring energy distance and the variogram metrics \cite{Pinson2012_EvaluatingQualityScenarios} between the ground truth (actual) wind power data and the scenarios.
    The project views support points as the critical link between forecasting and uncertainty quantification on the one hand, and the risk-aware optimizations on the other hand.

%% file: tex/optimization.tex
\section{OPTIMIZATION}
\label{sec:optimization}

As the penetration of renewable generation and distributed resources increases, so does the operational uncertainty that ISOs face, thereby calling for new operational paradigms.
In that context, this section discusses MISO's optimization pipeline, from day-ahead to real-time operations, under the combined lens of optimization, uncertainty, and risk.

Section \ref{sec:optimization:deterministic} describes current operational practices.
Section \ref{sec:optimization:stochastic} presents the project's plans for risk-aware optimization.
Section \ref{sec:optimization:extended} discusses multi-period market extensions.

\subsection{Current (deterministic) operations}
\label{sec:optimization:deterministic}

    \begin{figure}
        \centering
        \includegraphics[width=0.95\columnwidth]{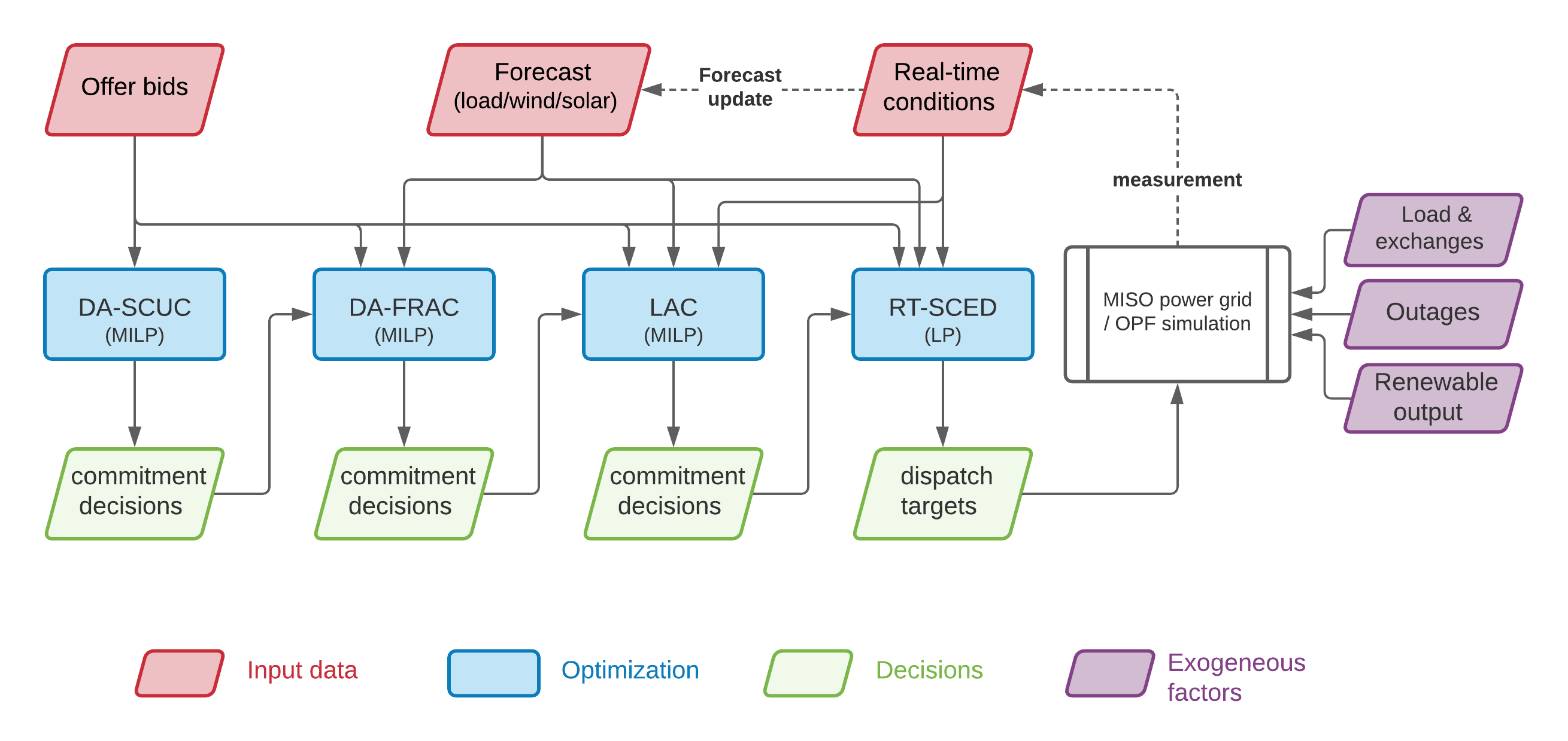}
        \caption{The MISO Optimization Pipeline, from Day-ahead to Real-Time.}
        \label{fig:optimization:miso_pipeline}
    \end{figure}

    From long-term planning and capacity expansion to clearing the real-time market, ISOs heavily rely on solving optimization problems.
    MISO's current operating pipeline, from day-ahead to real-time operations, is shown in Figure \ref{fig:optimization:miso_pipeline} and fully implemented (as a digital twin) as part of RAMC.

\begin{figure}[!t]
        \centering
        \includegraphics[width=0.9\columnwidth]{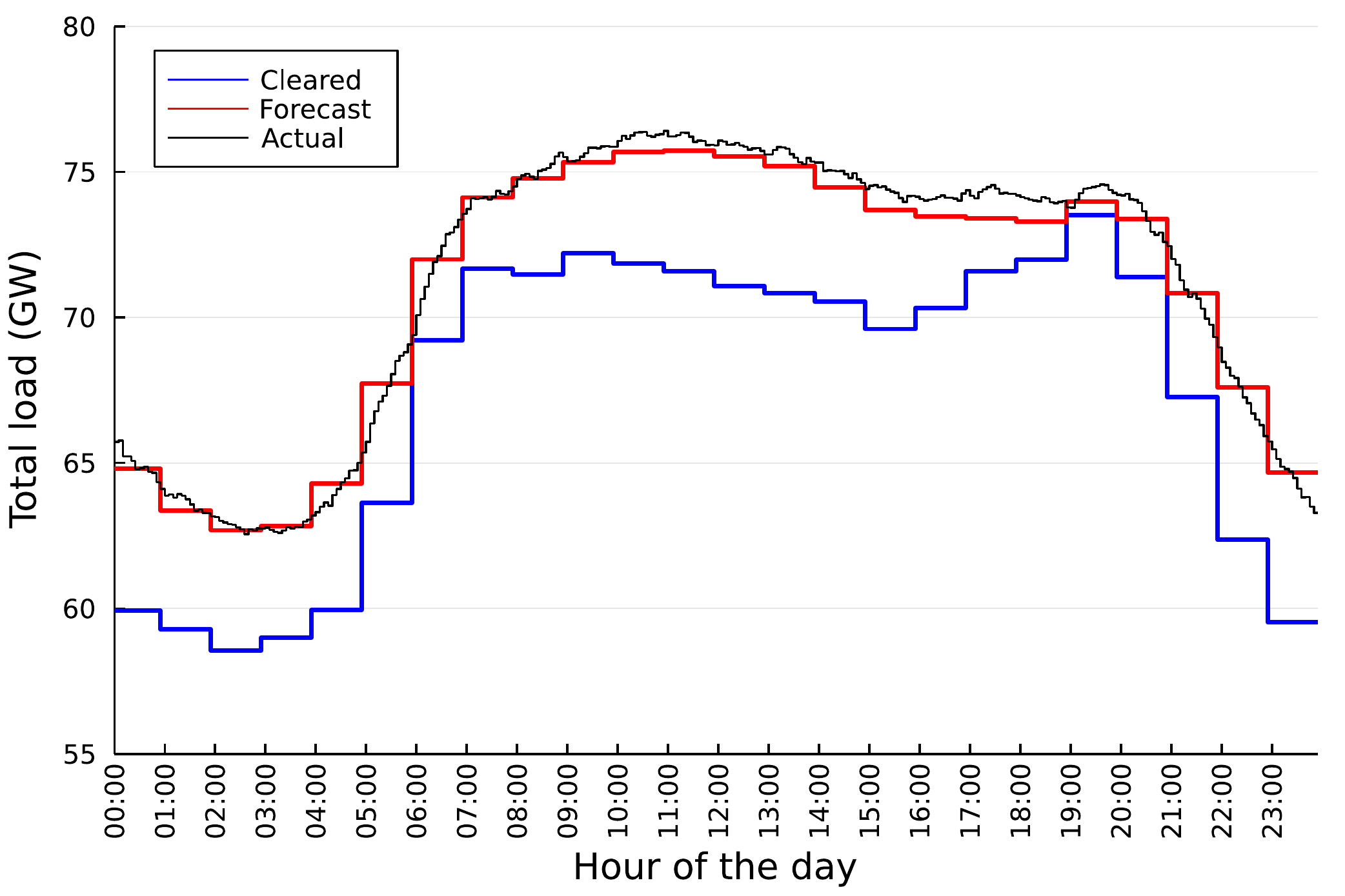}
        \caption{Total load in MISO throughout March 30th, 2022.
        \textcolor{blue}{\textbf{Cleared}}: total demand cleared by the DA-SCUC (hourly).
        \textcolor{red}{\textbf{Forecast}}: day-ahead load forecast used in DA-FRAC (hourly).
        \textcolor{black}{\textbf{Actual}}: actual load measured by MISO (5 minute).
        }
        \label{fig:optimization:MISO_load}
    \end{figure}

    The operations are organized as follows.
    First, between 10am and 1pm prior to the operating day, the day-ahead (DA) energy and reserve market is cleared.
    The DA market (DA-SCUC) takes supply and demand bids as inputs from each market participant, and produces commitment decisions for each generator.
    Second, the Forward Reliability Assessment Commitment (DA-FRAC) takes, as input, the previous commitment decisions and current forecasts for load and renewable outputs, and may commit additional generation units to serve the forecasted demand reliably.
    Figure \ref{fig:optimization:MISO_load} illustrates the difference between cleared, forecasted, and actual demand in MISO, on March 30\textsuperscript{th}, 2022.
    The demand cleared by DA-SCUC underestimates the forecast by 4\% on average, and by as much as 5GW; relative differences of $\pm 2{-}4\%$ are not uncommon throughout the year.
    Such discrepancies are addressed by committing additional units in DA-FRAC.

    Third, throughout the operating day, a Look-Ahead Commitment (LAC) optimization is executed, which again may commit additional units to accommodate previously unforeseen changes in net load.
    In MISO, LAC studies are executed every 15 minutes on a rolling-horizon fashion, with a 3-hour look-ahead window.
    Fourth, every 5 minutes, MISO's real-time market is cleared by solving a Security-Constrained Economic Dispatch problem (RT-SCED), which computes active power dispatch decisions for all (online) generators.
    Finally, dispatch decisions are sent to (human) operators who, after performing additional AC power flow studies, may adjust dispatch targets as appropriate.
    
    The DA-SCUC, DA-FRAC and LAC are security-constrained unit-commitment (SCUC) problems, and are modeled as mixed-integer linear programming (MILP) problems whose formulations are described in \cite{BPM_002}. The RT-SCED is a linear programming (LP) problem, based on a DC approximation with line losses \cite{BPM_002,Ma2009_MISOSCED}.
    Note that DA-SCUC, DA-FRAC and LAC only output commitment decisions, i.e., dispatch decisions are not final until the RT-SCED.
    The DA-SCUC, DA-FRAC, LAC and RT-SCED problems must be solved within 20, 40, 10 and 2 minutes, respectively.
    Despite significant progress in the resolution of UC problems \cite{Knueven2019_MILPFormulationsforSCUC}, industrial cases remain computationally and numerically challenging to solve \cite{Anjos2017_UCinElectricEnergySystems,Chen2021_HIPPO}.

\subsection{Future (stochastic) operations}
\label{sec:optimization:stochastic}

    A core component of the RAMC project is the use of stochastic formulations within reliability assessment studies. Indeed, while the DA-SCUC and RT-SCED problems are intrinsically deterministic problems used for market-clearing, reliability studies such as DA-FRAC and LAC seek to address uncertainty in future load and renewable outputs.
    
    The DA-FRAC is executed in a day-ahead fashion, and the commitment status of slow-start units cannot be modified on short notice throughout the operating day.
    This naturally yields a two-stage stochastic programming (TSSP) formulation, where the first- and second-stage decision variables are the (binary) commitment and (continuous) dispatch decisions, respectively; see also \cite{Van2018_LargeScaleUCUncertainty}.
    The project has implemented a TSSP formulation, solved by Bender's decomposition, where the uncertainty in load and renewables is captured via a finite number of scenarios with support points-based down-sampling.
    
\begin{figure}[!t]
        \centering
        \includegraphics[width=0.48\columnwidth]{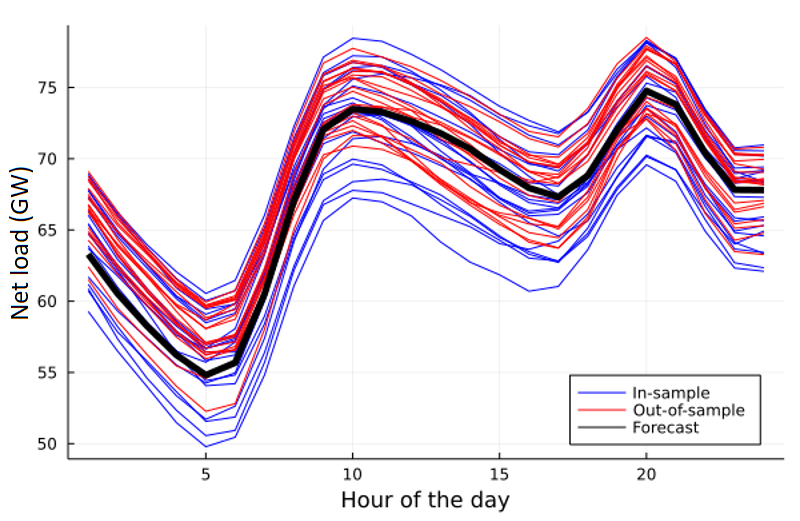}
        \hfill
        \includegraphics[width=0.48\columnwidth]{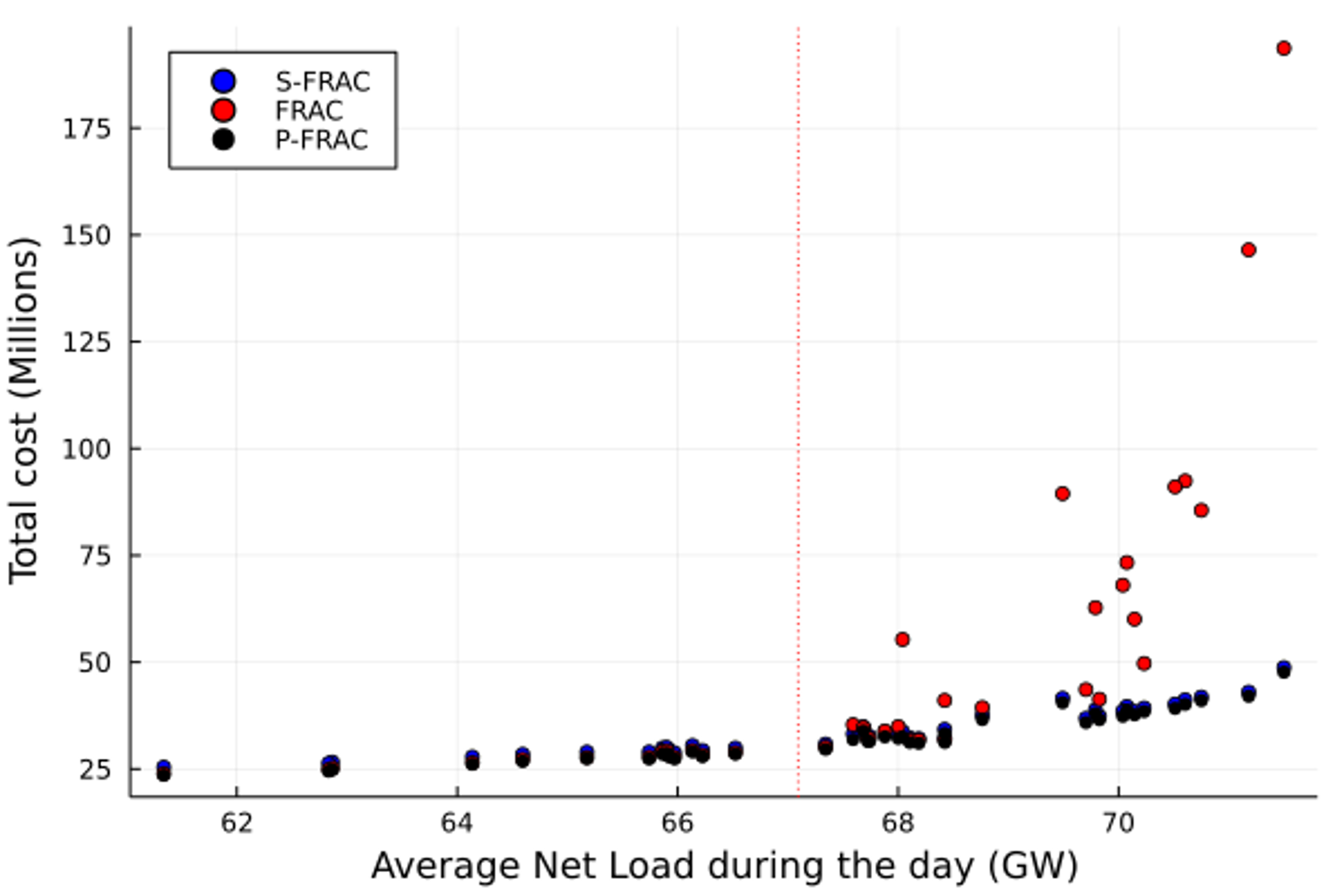}
        \caption{S-FRAC study in the RTE system on Feb. 2\textsuperscript{nd}, 2018.
        Left: net load forecast and scenarios throughout the day.
        Right: total operating costs for each scenario. Scenarios are ordered by average net load; the vertical line depicts the forecast's average net load.}
        \label{fig:optimization:SFRAC}
    \end{figure}

    The benefits of a stochastic FRAC have been evaluated on the French system, throughout an entire year.
    The results are illustrated in Figure \ref{fig:optimization:SFRAC}, which depicts net-load scenarios and compares total operating costs between deterministic (FRAC), stochastic (S-FRAC), and perfect-information (P-FRAC) formulations.
    The latter formulation relaxes non-anticipatory constraints, and allows to evaluate the value of perfect information (VPI).
    Results show that S-FRAC yields significant savings in high-net-load scenarios, for which S-FRAC achieves fewer and lower transmission violations.
    These reliability benefits result in monetary savings, as the ISO avoids costly corrective actions.
    Over the year, the VPI is reduced from 11\% in FRAC to 3\% in S-FRAC.
    
    LAC studies are executed in a rolling-horizon fashion, and potentially commit additional units throughout the day, i.e., current commitment decisions may be modified in the future.
    This process can be approximated as a TSSP problem wherein first-stage decisions are defined according to each resource's notification time, see, e.g., \cite{Knueven2021_SLAC}.
    The main difference with S-FRAC is that second-stage decisions may include discrete variables, which requires special treatment.
    State-of-the-art implementations typically leverage scenario-decomposition methods such as progressive hedging \cite{Van2018_LargeScaleUCUncertainty,Knueven2021_SLAC,Knueven2021_SOTAforStochasticUC}.
    Preliminary results on real-size instances from MISO suggest that stochastic LAC instances can be solved efficiently, echoing the findings of \cite{Knueven2021_SLAC}.

\subsection{Multi-Period Markets}
\label{sec:optimization:extended}
    
    In addition to stochastic optimization-based reliability studies, the project is investigating the benefits of a multi-period real-time market, also known as Look-Ahead Dispatch (LAD).
    Indeed, MISO's current RT-SCED is a single-period economic dispatch, which can lead to sub-optimal decisions when considering hours-long operations.
    LAD-based markets are also being investigated by other US ISOs  \cite{Xie2011_LookAheadDispatchERCOT,Hui2013_LookAheadUnforeseenERCOT,Hua2019_PricingMultiPeriodMarkets,Zhao2020_MultiPeriodMarket}.
    
    The project has implemented an hour-long LAD consisting of twelve five-minute time periods, and executed it on the French RTE system.
    Preliminary results demonstrate the computational feasibility of the approach, i.e., the LAD problems can be solved within the real-time market timeframe.
    Furthermore, the potential savings are significant: over the year, assuming a perfect forecast oracle for load and renewable, the LAD-based real-time market yields an overall cost reduction of 1.35\% compared to single-period SCED.
    Such savings result from lower transmission violations and lower operating costs, thanks to the better pre-positioning of generators' dispatch, especially during high-ramping hours.
    
    Nevertheless, additional experiments also suggest that large forecast errors can result in sub-optimal decisions in LAD compared to SCED.
    This effect increases with the length of the time horizon in LAD, and naturally calls for a stochastic LAD formulation, wherein uncertainty is explicitly taken into account.
    Another challenge in LAD is the accurate modeling of line losses in future time periods.
    The current RT-SCED formulation used by MISO relies on the so-called loss factors-based linearization \cite{Litvinov2004_LossLMP,FERC_LossApproximation}, which is only accurate within a small neighborhood of a given initial point.
    To overcome this limitation, the project has proposed a new linearization technique for line losses, which does not require any starting point \cite{Zhao2021_LLOA}.

%% file: tex/learning.tex
\section{LEARNING}
\label{sec:learning}

The increasing complexity of te market-clearing algorithms is accompanied by a growing need to quickly evaluate a system's behavior under a large numbers of scenarios for risk analysis. Therefore, the RAMC project is developing a new generation of tools that combine of ML and optimization techniques to address the above challenges. These include a novel \emph{just-in-time} learning framework, new architectures that exploit domain knowledge, and fast predict-then-repair procedures to compute close-to-optimal solutions.

\subsection{Just-In-Time Training}
\label{sec:learning:JIT}

    ML models that attempt to predict the solution of an Optimal Power Flow (OPF) problem have received a lot of attention recently \cite{Fioretto2020_PredictingACOPF,Chatzos2021_SpatialDecomposition}.
    Nevertheless, the vast majority of existing works only consider variability in demand and/or renewable output, and do not account for variability in commitment decisions nor generators' offers.

\begin{figure}[!t]
        \centering
        \includegraphics[width=0.9\columnwidth]{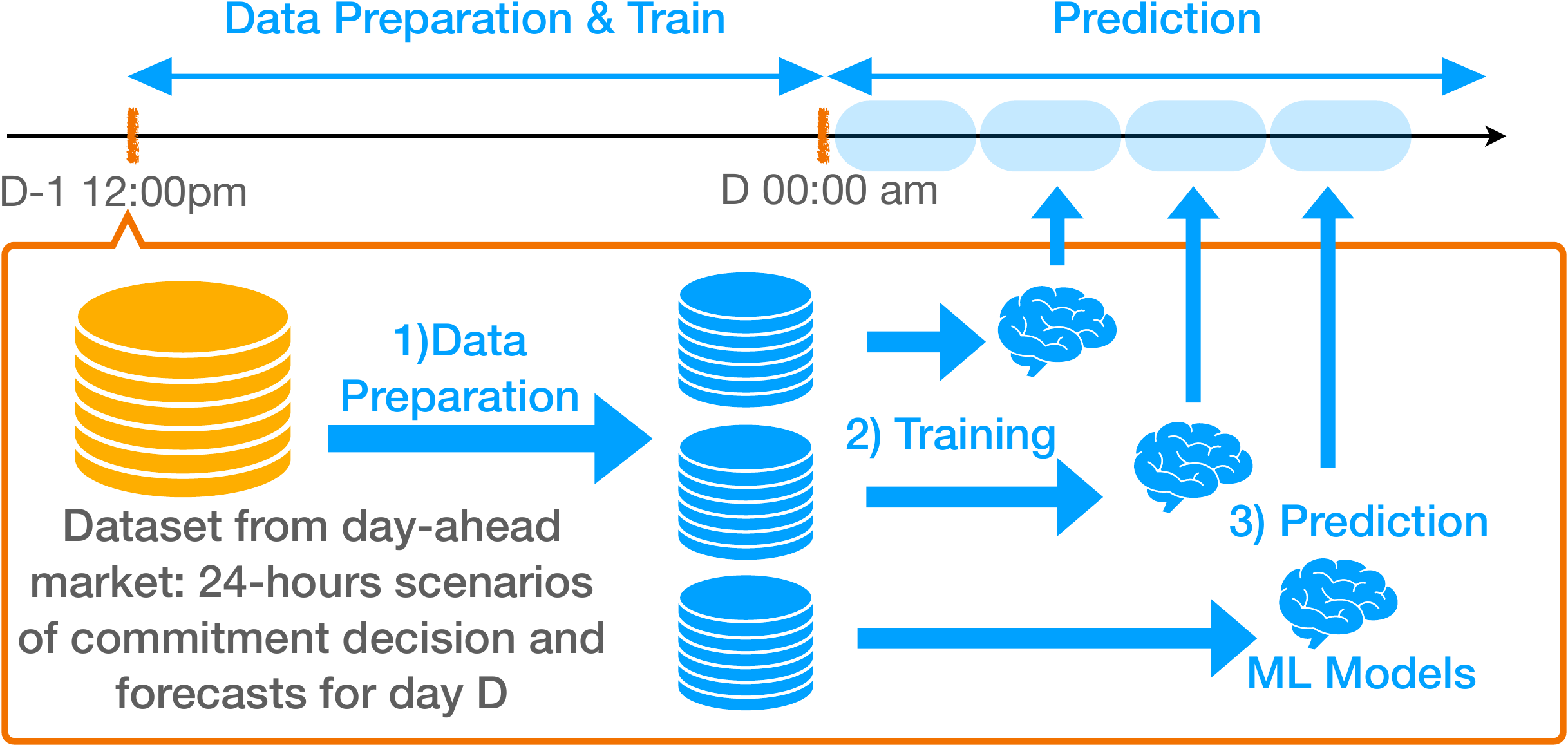}
        \caption{The proposed just-in-time machine learning framework. Excerpt from Fig. 1 in \cite{Chen2021_SCEDLearning}.}
        \label{fig:learning:JIT:pipeline}
    \end{figure}

    
    To tackle the combinatorial explosion of commitment decisions and the generator offers' variability, the project introduced in \cite{Chen2021_SCEDLearning} the novel concept of just-in-time training.
    The framework leverages the fact that, once the DA-SCUC and DA-FRAC are solved, $99\%$ of final commitment decisions are known \cite{Chen2018_MarketClearingSoftware}.
    In addition, day-ahead forecasting and scenario-generation tools provide an accurate distribution of future load and renewable output.
    The training process is illustrated in Figure \ref{fig:learning:JIT:pipeline}, and it proceeds as follows.
    
    First, hourly commitment decisions and generators' offers are gathered from the results of the DA-SCUC and DA-FRAC, and a large number of scenarios are generated for load and renewable outputs.
    Second, for each such scenario, a day-long simulation is executed, wherein an RT-SCED is solved every 5 minutes on a rolling basis.
    This provides a dataset of RT-SCED instances and corresponding solutions.
    Third, separate models, one for each hour of the day, are trained in a distributed fashion.
    Fourth, throughout the operating day, the trained models are used to make predictions.

    The main advantages of just-in-time training are twofold. On the one hand, it eliminates the combinatorial explosion in commitment decisions, because training takes place \emph{after} commitment decisions are made.
    On the other hand, it reduces the potential distribution shift between training and test data, because load and renewable scenarios are directly sampled from a day-ahead forecast model, as opposed to historical data.
    Finally, as reported in \cite{Chen2021_SCEDLearning}, the combined data-generation and training steps take no longer than a few hours, which is compatible with MISO's existing operations.

\subsection{Physics-Informed Architectures}
\label{sec:learning:architectures}
    
    The project has developed several methodologies to improve the accuracy and training efficiency of ML models, leveraging domain knowledge and the intrinsic structure found in power systems.
    A Lagrangian loss function has been introduced in \cite{Fioretto2020_PredictingACOPF}, to mitigate the fact that ML-based predictions may --and usually do-- violate physical constraints.
    The proposed Lagrangian loss penalizes constraint violation during training. The training is performed using a Lagrangian Dual approach that significantly improves the feasibility of the predicted AC-OPF solutions.
    This approach has been successfully employed in OPF settings \cite{Velloso2021_DNNSCOPF,Chatzos2021_SpatialDecomposition}.
    
    \begin{figure}[!t]
        \centering
        \includegraphics[width=\columnwidth]{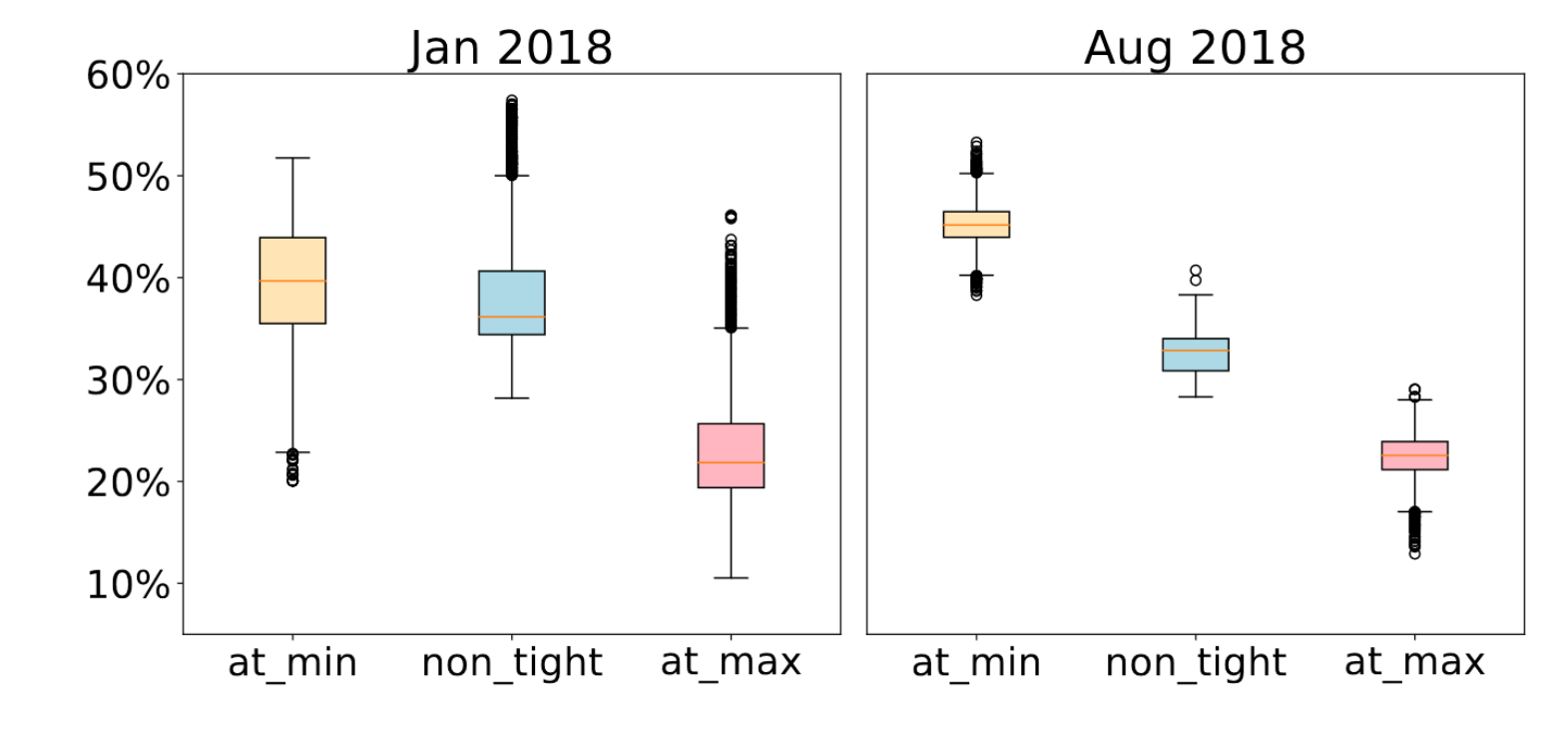}
        \caption{Proportion of non-renewable generators at their minimum/maximum limits in January (left) and August (right). Excerpt from Fig. 3 in \cite{Chen2021_SCEDLearning}.}
        \label{fig:learning:architectures:active_bounds}
    \end{figure}
    
    Novel learning architectures have also been proposed. In the context of learning AC-OPF solutions, a spatial decomposition-based architecture has been introduced in \cite{Chatzos2021_SpatialDecomposition}, in conjunction with a distributed, two-stage training procedure. Again, the innovation here is to transfer the optimization concept of Lagrangian decomposition to machine learning. 
    This architecture and training method directly exploit the spatial structure, namely, regional decomposition, typically found in power grids.
    Furthermore, in \cite{Chen2021_SCEDLearning}, is was observed that most generators are dispatched to either their minimum or maximum limit, as is illustrated in Figure \ref{fig:learning:architectures:active_bounds}.
    This motivated a classification-then-regression architecture, wherein a classifier first predicts whether lower/upper bound constraints are active, and a regression model predicts the dispatch of the remaining generators.
    Extensive numerical results demonstrate the benefits of the combined method over regression models alone.

\subsection{Predict Then Repair}
\label{sec:learning:ml4opt}

    A fundamental limitation of optimization proxies is the absence of guarantee that the predicted output satisfies the constraints of the optimization problem.
    For instance, a predicted dispatch may violate the power balance constraint, or a generator's predicted commitment schedule may not respect its minimum down time.
    Although an approximate solution is sometimes acceptable, e.g., when approximating total operating costs or checking whether large violations occur, accurate and feasible solutions are often required.
    To overcome this challenge, the project implemented fast feasibility-restoration procedures that generate close-to-optimal feasible solutions, when seeded with an almost-feasible and almost-optimal prediction.
    This latter assumption is crucial to  practical performance.
    
    A classical approach is to project the predicted solution onto the feasible set of the original problem. This technique is employed in \cite{Velloso2021_DNNSCOPF} and \cite{Chatzos2021_SpatialDecomposition}, in the context of optimal power flow problems.
    Nevertheless, despite strong guarantees, doing so may be as hard as solving the original problem, especially when dealing with non-convexities. Thus, in \cite{Chatzos2021_SpatialDecomposition}, the authors considered the faster alternative of solving a power flow problem to restore feasibility, and report an order of magnitude improvement over the projection approach.
    One may also predict partial solutions, or sub-spaces that are likely to contain an optimal solution.
    A reduced problem is then solved, hopefully significantly faster than the original one.
    Experiments on unit commitment problems are reported in, e.g., \cite{Mossina2019_MultiLabelClassification,Xavier2021_LearningLargeScaleSCUC}, and suggest that high-quality solutions can be obtained with a five- to ten-fold speedup.
    
    Finally, the above techniques can naturally be used as heuristics within a traditional optimization algorithm.
    Nevertheless, it is expected that such techniques are useful even outside of optimization frameworks, under strict time budget.
    The reader is referred to \cite{Bengio2021_ML4COMethodological} for a broad overview of the interplay between ML and optimization.

%% file: tex/risk.tex
\section{RISK-AWARE OPERATIONS}
\label{sec:risk}

\begin{figure}
    \centering
    \includegraphics[width=\columnwidth]{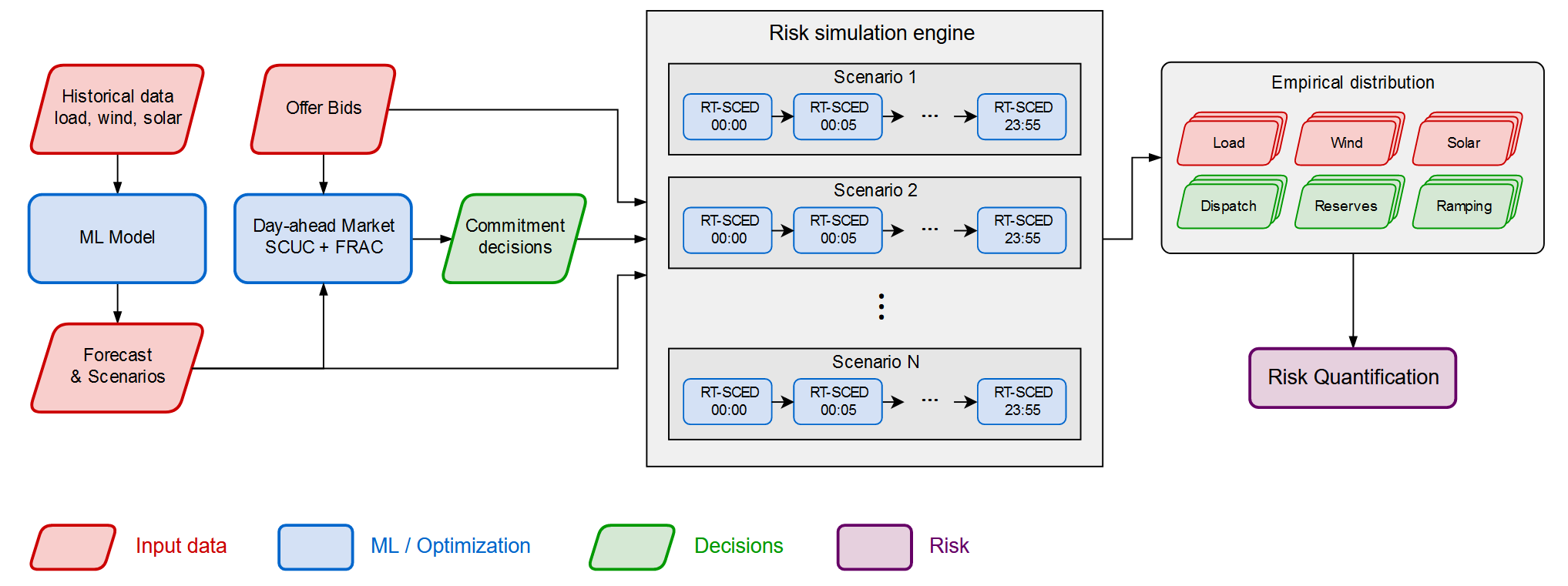}
    \caption{Risk estimation pipeline. For given day-ahead commitments, multiple day-long simulations are executed, which produce empirical distributions for load, renewable output, reserves, etc. These empirical distributions are then aggregated into risk metrics.}
    \label{fig:risk:pipeline}
\end{figure}

The combination of forecasting and uncertainty-quantification models, digital twins, and accurate ML-based optimization proxies, is the foundation for systematic risk quantification and risk-aware operations, which is the unifying goal of the RAMC project. First, the project is developing new risk metrics that quantify risk in a principle fashion. Such metrics include the conditional expectation of a quantity of interest, e.g., total ramping capacity, over the $\alpha\%$ worst cases, also known as \emph{conditional value at risk} (CVaR), the probability of an adverse event occurring, and the expected financial cost of such adverse events.
The risk estimation pipeline, as implemented in RAMC, is depicted in Figure \ref{fig:risk:pipeline}: simulation runs are evaluated in parallel, and the resulting empirical distributions are used to quantify risk. While Figure \ref{fig:risk:pipeline} illustrates a day-ahead setting, the risk pipeline can accommodate arbitrary study window, e.g., a few hours ahead, and recourse actions such as additional commitments from LAC.

Although the risk estimation pipeline leverages parallel computing, solving a large number of optimization problems remains costly, especially for real-time decision making.
For instance, MISO's real-time market is cleared every five minutes, which means that operators have a very short time-window, no more than a couple of minutes, to evaluate the proposed dispatch targets and make necessary adjustments. In that context, tools are needed that can provide accurate risk quantification within a second, thereby precluding the use optimization-based simulations.
To overcome this limitation, the project makes use of the aforementioned ML-based optimization proxies, which can evaluate hundreds of thousands of scenarios in a second \cite{Chen2021_SCEDLearning}.
The ML-based simulations are then aggregated into risk metrics, which allows operator feedback in real-time.

Finally, risk metrics such as CVaR can be incorporated in the ISO's optimization pipeline, resulting in risk-aware market-clearing algorithms. Doing so remains computationally tractable, because CVaR can be represented as a convex function. The resulting formulations enable the exploration of the trade-off between risk and cost, and ensure reliable operations without significant cost increases.

%% file: tex/conclusion.tex
\section{CONCLUSION}
\label{sec:conclusion}

This paper presented several challenges that arise as electrical power grids transition from fossil fuel-based to renewable energy generation sources. The combined increase in demand stochasticity and production volatility calls for a new operational paradigm, that explicitly quantifies risk and embeds it in the algorithms that drive the power grid. The paper described the unified risk-aware market-clearing framework pursued in the RAMC project. In particular, the project has investigated new forecasting and uncertainty-quantification tools, advocating that support points provide a unique bridge between forecasting and risk-aware optimization. Furthermore, the computational feasibility of embedding stochastic optimization techniques within an ISO's optimization pipeline has been demonstrated on real, industrial-size systems. Experimental results also suggest significant benefits can be achieved by doing so. New learning frameworks have been also proposed, from just-in-time training to the Lagrangian dual training of physics-informed architectures and a novel "Predict then Repait" framework, demonstrating that market-clearing algorithms can significant benefits from data-driven approach to speed up computations. Finally, these advances in prediction, optimization and learning have been integrated into a unified risk-estimation pipeline through novel methodologies. Overall, the combination of new risk metrics, stochastic optimization tools, and fast ML-based optimization proxies, represents a promising direction towards a principled risk-aware operational framework for ISOs.

%% file: root.bbl
\begin{thebibliography}{10}
\providecommand{\url}[1]{#1}
\csname url@rmstyle\endcsname
\providecommand{\newblock}{\relax}
\providecommand{\bibinfo}[2]{#2}
\providecommand\BIBentrySTDinterwordspacing{\spaceskip=0pt\relax}
\providecommand\BIBentryALTinterwordstretchfactor{4}
\providecommand\BIBentryALTinterwordspacing{\spaceskip=\fontdimen2\font plus
\BIBentryALTinterwordstretchfactor\fontdimen3\font minus
  \fontdimen4\font\relax}
\providecommand\BIBforeignlanguage[2]{{%
\expandafter\ifx\csname l@#1\endcsname\relax
\typeout{** WARNING: IEEEtran.bst: No hyphenation pattern has been}%
\typeout{** loaded for the language `#1'. Using the pattern for}%
\typeout{** the default language instead.}%
\else
\language=\csname l@#1\endcsname
\fi
#2}}

\bibitem{EIA860}
\BIBentryALTinterwordspacing
{US Energy Information Administration}, ``{Form EIA-860},'' 2020. [Online].
  Available: \url{https://www.eia.gov/electricity/data/eia860/}
\BIBentrySTDinterwordspacing

\bibitem{ODRE}
\BIBentryALTinterwordspacing
{Open Data Réseaux Énergies}, ``{Registre national des installations de
  production et de stockage d'électricité 2018}.'' [Online]. Available:
  \url{https://opendata.reseaux-energies.fr/}
\BIBentrySTDinterwordspacing

\bibitem{Zhu2021_MultiResolutionSpatioTemporal}
S.~Zhu, H.~Zhang, Y.~Xie, and P.~Van~Hentenryck, ``Multi-resolution
  spatio-temporal prediction with application to wind power generation,''
  \emph{arXiv preprint arXiv:2108.13285}, 2021.

\bibitem{ursu2009modelling}
E.~Ursu and P.~Duchesne, ``On modelling and diagnostic checking of vector
  periodic autoregressive time series models,'' \emph{Journal of Time Series
  Analysis}, vol.~30, no.~1, pp. 70--96, 2009.

\bibitem{hu2018nonlinear}
Y.-L. Hu and L.~Chen, ``A nonlinear hybrid wind speed forecasting model using
  lstm network, hysteretic elm and differential evolution algorithm,''
  \emph{Energy conversion and management}, vol. 173, pp. 123--142, 2018.

\bibitem{chung2014empirical}
J.~Chung, C.~Gulcehre, K.~Cho, and Y.~Bengio, ``Empirical evaluation of gated
  recurrent neural networks on sequence modeling,'' \emph{arXiv preprint
  arXiv:1412.3555}, 2014.

\bibitem{sweeney2020future}
C.~Sweeney, R.~J. Bessa, J.~Browell, and P.~Pinson, ``The future of forecasting
  for renewable energy,'' \emph{Wiley Interdisciplinary Reviews: Energy and
  Environment}, vol.~9, no.~2, p. e365, 2020.

\bibitem{yu2020superposition}
M.~Yu, Z.~Zhang, X.~Li, J.~Yu, J.~Gao, Z.~Liu, B.~You, X.~Zheng, and R.~Yu,
  ``Superposition graph neural network for offshore wind power prediction,''
  \emph{Future Generation Computer Systems}, vol. 113, pp. 145--157, 2020.

\bibitem{park2019physics}
J.~Park and J.~Park, ``Physics-induced graph neural network: An application to
  wind-farm power estimation,'' \emph{Energy}, vol. 187, p. 115883, 2019.

\bibitem{chatzos2021data}
M.~Chatzos, M.~Tanneau, and P.~Van~Hentenryck, ``Data-driven time series
  reconstruction for modern power systems research,'' \emph{arXiv preprint
  arXiv:2110.13772}, 2021.

\bibitem{gal2016dropout}
Y.~Gal and Z.~Ghahramani, ``Dropout as a bayesian approximation: Representing
  model uncertainty in deep learning,'' in \emph{international conference on
  machine learning}.\hskip 1em plus 0.5em minus 0.4em\relax PMLR, 2016, pp.
  1050--1059.

\bibitem{werho2021scenario}
T.~Werho, J.~Zhang, V.~Vittal, Y.~Chen, A.~Thatte, and L.~Zhao, ``Scenario
  generation of wind farm power for real-time system operation,'' \emph{arXiv
  preprint arXiv:2106.09105}, 2021.

\bibitem{zhang2019short}
J.~Zhang, J.~Yan, D.~Infield, Y.~Liu, and F.-s. Lien, ``Short-term forecasting
  and uncertainty analysis of wind turbine power based on long short-term
  memory network and gaussian mixture model,'' \emph{Applied Energy}, vol. 241,
  pp. 229--244, 2019.

\bibitem{dumas2022deep}
J.~Dumas, A.~Wehenkel, D.~Lanaspeze, B.~Corn{\'e}lusse, and A.~Sutera, ``A deep
  generative model for probabilistic energy forecasting in power systems:
  normalizing flows,'' \emph{Applied Energy}, vol. 305, p. 117871, 2022.

\bibitem{zheng2014stochastic}
Q.~P. Zheng, J.~Wang, and A.~L. Liu, ``Stochastic optimization for unit
  commitment—a review,'' \emph{IEEE Transactions on Power Systems}, vol.~30,
  no.~4, pp. 1913--1924, 2014.

\bibitem{mak2018support}
S.~Mak and V.~R. Joseph, ``Support points,'' \emph{The Annals of Statistics},
  vol.~46, no.~6A, pp. 2562--2592, 2018.

\bibitem{Vakayil22_DataTwinning}
A.~Vakayil and V.~R. Joseph, ``Data twinning,'' \emph{Statistical Analysis and
  Data Mining: The ASA Data Science Journal}, 2022.

\bibitem{chen2013wind}
N.~Chen, Z.~Qian, I.~T. Nabney, and X.~Meng, ``Wind power forecasts using
  gaussian processes and numerical weather prediction,'' \emph{IEEE
  Transactions on Power Systems}, vol.~29, no.~2, pp. 656--665, 2013.

\bibitem{Pinson2012_EvaluatingQualityScenarios}
P.~Pinson and R.~Girard, ``Evaluating the quality of scenarios of short-term
  wind power generation,'' \emph{Applied Energy}, vol.~96, pp. 12--20, 2012,
  {S}mart {G}rids.

\bibitem{BPM_002}
MISO, ``Energy and operating reserve markets,'' {B}usiness Practices Manual
  Energy and Operating Reserve Markets.

\bibitem{Ma2009_MISOSCED}
X.~Ma, Y.~Chen, and J.~Wan, ``Midwest {ISO} co-optimization based real-time
  dispatch and pricing of energy and ancillary services,'' in \emph{2009 IEEE
  Power Energy Society General Meeting}, 2009, pp. 1--6.

\bibitem{Knueven2019_MILPFormulationsforSCUC}
B.~Knueven, J.~Ostrowski, and J.-P. Watson, ``On mixed-integer programming
  formulations for the unit commitment problem,'' \emph{INFORMS Journal on
  Computing}, vol.~32, no.~4, pp. 857--876, 2020.

\bibitem{Anjos2017_UCinElectricEnergySystems}
M.~F. Anjos, A.~J. Conejo, \emph{et~al.}, \emph{Unit commitment in electric
  energy systems}.\hskip 1em plus 0.5em minus 0.4em\relax Now Publishers, Inc.,
  2017, vol.~1, no.~4.

\bibitem{Chen2021_HIPPO}
Y.~Chen, F.~Pan, J.~Holzer, E.~Rothberg, Y.~Ma, and A.~Veeramany, ``A high
  performance computing based market economics driven neighborhood search and
  polishing algorithm for security constrained unit commitment,'' \emph{IEEE
  Transactions on Power Systems}, vol.~36, no.~1, pp. 292--302, 2021.

\bibitem{Van2018_LargeScaleUCUncertainty}
W.~van Ackooij, I.~Danti~Lopez, A.~Frangioni, F.~Lacalandra, and M.~Tahanan,
  ``Large-scale unit commitment under uncertainty: an updated literature
  survey,'' \emph{Annals of Operations Research}, vol. 271, no.~1, pp. 11--85,
  2018.

\bibitem{Knueven2021_SLAC}
B.~Knueven, M.~Faqiry, M.~Garcia, Y.~Chen, R.~Treinen, T.~Werho, J.~Zhang,
  V.~Vittal, L.~Zhao, A.~Thatte, \emph{et~al.}, ``Stochastic look-ahead
  commitment: A case study in miso,'' 2021.

\bibitem{Knueven2021_SOTAforStochasticUC}
\BIBentryALTinterwordspacing
B.~Knueven, J.-P. Watson, J.~Ostrowski, and D.~L. Woodruff, ``State-of-the-art
  techniques for large-scale stochastic unit commitment,'' 6 2021. [Online].
  Available: \url{https://www.osti.gov/biblio/1804443}
\BIBentrySTDinterwordspacing

\bibitem{Xie2011_LookAheadDispatchERCOT}
L.~Xie, X.~Luo, and O.~Obadina, ``Look-ahead dispatch in ercot: Case study,''
  in \emph{2011 IEEE Power and Energy Society General Meeting}, 2011, pp. 1--3.

\bibitem{Hui2013_LookAheadUnforeseenERCOT}
H.~Hui, C.-N. Yu, R.~Surendran, F.~Gao, S.~Moorty, and X.~Xu, ``Look ahead to
  the unforeseen: Ercot's nonbinding look-ahead sced study,'' in \emph{2013
  IEEE Power Energy Society General Meeting}, 2013, pp. 1--5.

\bibitem{Hua2019_PricingMultiPeriodMarkets}
B.~Hua, D.~A. Schiro, T.~Zheng, R.~Baldick, and E.~Litvinov, ``Pricing in
  multi-interval real-time markets,'' \emph{IEEE Transactions on Power
  Systems}, vol.~34, no.~4, pp. 2696--2705, 2019.

\bibitem{Zhao2020_MultiPeriodMarket}
J.~Zhao, T.~Zheng, and E.~Litvinov, ``A multi-period market design for markets
  with intertemporal constraints,'' \emph{IEEE Transactions on Power Systems},
  vol.~35, no.~4, pp. 3015--3025, 2020.

\bibitem{Litvinov2004_LossLMP}
E.~Litvinov, T.~Zheng, G.~Rosenwald, and P.~Shamsollahi, ``Marginal loss
  modeling in lmp calculation,'' \emph{IEEE Transactions on Power Systems},
  vol.~19, no.~2, pp. 880--888, 2004.

\bibitem{FERC_LossApproximation}
B.~Eldridge, R.~P. O’Neill, and A.~Castillo, ``Marginal loss calculations for
  the dcopf,'' 2017.

\bibitem{Zhao2021_LLOA}
H.~Zhao, M.~Tanneau, and P.~Van~Hentenryck, ``{A Linear Outer Approximation of
  Line Losses for DC-based Optimal Power Flow Problems},'' 2021.

\bibitem{Fioretto2020_PredictingACOPF}
F.~Fioretto, T.~W. Mak, and P.~Van~Hentenryck, ``Predicting ac optimal power
  flows: Combining deep learning and lagrangian dual methods,'' in
  \emph{Proceedings of the AAAI Conference on Artificial Intelligence},
  vol.~34, no.~01, 2020, pp. 630--637.

\bibitem{Chatzos2021_SpatialDecomposition}
M.~Chatzos, T.~W. Mak, and P.~Vanhentenryck, ``Spatial network decomposition
  for fast and scalable ac-opf learning,'' \emph{IEEE Transactions on Power
  Systems}, pp. 1--1, 2021.

\bibitem{Chen2021_SCEDLearning}
W.~Chen, S.~Park, M.~Tanneau, and P.~Van~Hentenryck, ``Learning optimization
  proxies for large-scale security-constrained economic dispatch,'' 2021.

\bibitem{Chen2018_MarketClearingSoftware}
Y.~Chen, F.~Wang, J.~Wan, and F.~Pan, ``Developing next generation electricity
  market clearing optimization software,'' in \emph{2018 IEEE Power Energy
  Society General Meeting (PESGM)}, 2018, pp. 1--5.

\bibitem{Velloso2021_DNNSCOPF}
A.~Velloso and P.~Van~Hentenryck, ``Combining deep learning and optimization
  for preventive security-constrained dc optimal power flow,'' \emph{IEEE
  Transactions on Power Systems}, vol.~36, no.~4, pp. 3618--3628, 2021.

\bibitem{Mossina2019_MultiLabelClassification}
L.~Mossina, E.~Rachelson, and D.~Delahaye, ``Multi-label classification for the
  generation of sub-problems in time-constrained combinatorial optimization,''
  2019.

\bibitem{Xavier2021_LearningLargeScaleSCUC}
A.~S. Xavier, F.~Qiu, and S.~Ahmed, ``Learning to solve large-scale
  security-constrained unit commitment problems,'' \emph{INFORMS Journal on
  Computing}, vol.~33, no.~2, pp. 739--756, 2021.

\bibitem{Bengio2021_ML4COMethodological}
Y.~Bengio, A.~Lodi, and A.~Prouvost, ``Machine learning for combinatorial
  optimization: a methodological tour d’horizon,'' \emph{European Journal of
  Operational Research}, vol. 290, no.~2, pp. 405--421, 2021.

\end{thebibliography}
